\documentclass{daj}

\usepackage{amsthm, amssymb, enumerate,amsmath}

\theoremstyle{plain}
\newtheorem{theorem}{Theorem}

\newtheorem{lemma}[theorem]{Lemma}
\newtheorem{example}[theorem]{Example}

\dajAUTHORdetails{%
  title = {A Question of Bukh on Sums of Dilates}, 
  author = {Brandon Hanson and Giorgis Petridis},
  plaintextauthor = {Brandon Hanson, Giorgis Petridis},
  keywords = {sumsets, dilate sums},
}   

\dajEDITORdetails{%
   year={2021},
   number={13},
   received={28 September 2020},   
   published={15 September 2021},  
   doi={10.19086/da.28143},       
}   

\begin{document}

\begin{frontmatter}[classification=text]


\author[han]{Brandon Hanson\thanks{Supported by the NSF Award 2001622.}}
\author[pet]{Giorgis Petridis\thanks{Supported by the NSF Award 1723016 and 2054214; and gratefully acknowledges the support from the RTG in Algebraic Geometry, Algebra, and Number Theory at the University of Georgia.}}

\begin{abstract}
We answer in the affirmative a question of Bukh on the cardinality of the dilate sum $A + 2 \cdot A$.
\end{abstract}
\end{frontmatter}


\section{Introduction}

For finite set of integers $A$ and a positive integer $\lambda$ define
\[
A+ \lambda \cdot A =\{ a + \lambda a' : a,a' \in A \}.
\]
Such dilate sums were investigated by Bukh in \cite{Bu2}, where the asymptotically sharp lower bound $|A+ \lambda \cdot A| \geq (\lambda+1) |A| - o(|A|)$ was obtained. The term $o(|A|)$ was replaced by $O_\lambda(1)$ in \cite{BaSh,Sha}. For certain $\lambda$ there are sharp lower bounds for this type of dilate sum \cite{CHS,CSV,DCS,HaRu,Na}. The question has also been investigated for finite sets of real numbers \cite{BrGr,ChFa} and of other commutative groups \cite{F-P,Mu, Pla,PlTr}.  

The trivial upper bound $|A + \lambda \cdot A| \leq |A|^2$ is attained when $\lambda=2$ and $A$ is a geometric progression with common ratio 3. However, both from a theoretical standpoint and in applications \cite{Bu1,Bu2,CrSi}, one is interested in bounding dilate sums in terms of $|A|$, $\lambda$, and the doubling parameter
\[
K = \frac{|A+A|}{|A|}.
\]

An inequality of Pl\"unnecke (see the remark after Lemma~\ref{Plunnecke}), coupled with the inclusion of $A +\lambda \cdot A$ in the $(\lambda+1)$-fold sumset of $A$, implies $|A +\lambda \cdot A| \leq K^{\lambda+1} |A|.$

Bukh improved this to $|A+ \lambda \cdot A| \leq K^{C \log\lambda} |A|$ \cite{Bu2}. Bukh's bound is qualitatively optimal, as we see by setting $A = \{1, \dots, n\}$ and $\lambda=n$. However, the method in \cite{Bu2}, and its refinement in \cite{BuZh}, do not improve the Pl\"unnecke-induced upper bound for small values of $\lambda$. 

The case $\lambda=2$ is illustrative. As we have seen the state-of-the-art follows from two inequalities:
\begin{equation}\label{Plu A+2A}
|A + 2 \cdot A| \leq |A+A+A| \leq K^3 |A|.    
\end{equation}
The first inequality above is sharp for generalised arithmetic progressions and is sharp up to multiplicative constants for dense random subsets of generalised arithmetic progressions. The second inequality is sharp up to multiplicative constants for some examples of Ruzsa \cite[Theorem 9.5]{Ru4}. The examples that show near sharpness for one of the two inequalities are far from sharp for the other inequality. This raises the possibility that \eqref{Plu A+2A} can be improved. 

Bukh asked \cite[Question 4.3]{Bu2} if there exists a $p<3$ such that for all $K \in \mathbb{R}$ and all finite $A \subset \mathbb{Z}$ such that $|A+A| \leq K |A|$, we have
\[
|A + 2 \cdot A | \leq K^p |A|?
\]

The main result of the present paper is an affirmative answer to Bukh's question. We show that one may take $p = 3-1/20$. Our result holds in any commutative group, where for two subsets $A$ and $B$ we define
\[
A+ 2 \cdot B = \{a+b+b : a \in A , b \in B\}.
\]

\begin{theorem}\label{Dilates}
Let $A$ be a subset of a commutative group $G$ and $K$ be a parameter. If $|A+A| \leq K |A|$, then $|A+2 \cdot A| \leq K^{2.95} |A|.$
\end{theorem}

The methods employed in the proof of Theorem~\ref{Dilates} combine widely used tools of additive combinatorics, see Section~\ref{Lemmas}, with techniques developed in the study of Fre\u{\i}man's theorem \cite{Fr}. More specifically work of Katz and Koester \cite{KaKo}, Schoen \cite{Sch}, and Schoen and Shkredov~\cite{ScSh}.

Taking $A = \{0,1\}^n \subset \mathbb{Z}^n$ (or a Fre\u{\i}man 3-isomorphic subset in $\mathbb{Z}$) shows that the exponent $p$ in Bukh's question must be at least $ \log2 / \log(3/2) > 1.7$.

\subsection*{Further results and organisation of the paper}

Some of the tools used in the proof of Theorem~\ref{Dilates} are listed in Section \ref{Lemmas}. The proof of Theorem~\ref{Dilates} is carried out in Sections \ref{Proof 1} and \ref{Proof2}. An overview can be found below Lemma~\ref{BasicLemma}.

Theorem~\ref{Dilates} can be extended to $A-2 \cdot A$ under the stronger hypothesis that both $|A+A|$ and $|A-A|$ are bounded above by $K|A|$, as can be assumed in the applications in \cite{Bu1,Bu2,CrSi}. This topic is discussed in Section~\ref{A-2A}.

The proof of Theorem~\ref{Dilates} works for all positive integer $\lambda >0$. It gives, to the best knowledge of the authors, new upper bounds on dilate sums $A+ \lambda A$ for $\lambda=2,3,5$. This topic is investigated in Section~\ref{A+ prime A}.

In Section~\ref{Final Section} we offer an improvement to Theorem~\ref{Dilates} for large $K$. Modifying arguments of Ruzsa from \cite{Ru3,Ru4}, we prove 
\begin{equation}\label{Large K}
|A+2 \cdot A| \leq (K |A|)^{4/3}.
\end{equation}
This is the best known bound when $|A|^{20/97} \leq K \leq |A|^{1/2}$.

\section{Lemmas} \label{Lemmas}

\subsection{Pl\"unnecke's inequality}

We make heavy use of an inequality of Pl\"unnecke that bounds the cardinality of higher sumsets in term of the doubling parameter \cite{Plu,Ru2}. We need the form proved in \cite{Pe}.

\begin{lemma}[Pl\"unnecke's inequality]\label{Plunnecke}
Let $U$ and $V$ be a finite subsets of a commutative group $G$ and let $X\subseteq U$ be such that \[\frac{|X+V|}{|X|}\leq \frac{|X'+V|}{|X'|}\] for all non-empty subsets $\emptyset \neq X'\subseteq X$. Then for any set $W$, we have
\[|X+V+W|\leq \frac{|X+V||X+W|}{|X|}.\]
In particular, the stated bound on $|X+V+W|$ holds if  
\[\frac{|X+V|}{|X|}\leq \frac{|X'+V|}{|X'|}\]
for all non-empty subsets $\emptyset \neq X'\subseteq U$ (because these include all subsets of $X$) .
\end{lemma}

By iterating the above we get that if $|A+A| \leq K |A|$, then 
\[
|A+A+A| \leq |X+A+A+A| \leq K^3 |X| \leq K^3 |A|.
\]

\subsection{Popular differences and the Katz--Koester inclusion}

We introduce notation that will be used throughout the paper. Given two finite subsets $U$ and $V$ of a commutative group $G$ and $d \in U-V$ we will often denote by $U_d$ the set $U \cap (d +V)$ (suppressing $V$). We will also abuse notation and say that $d$ is a $t$-popular difference in $U - V$ when $|U_d| \geq t$.

A simple fact we will use repeatedly is
\[
\sum_{d \in U - V} |U \cap (d+V)| = |U| |V|.
\]

Katz and Koester showed the importance of (a variant of) the following simple observation in \cite{KaKo}. Given any finite subset $W \subseteq G$ and $d \in U-V$ 
\[
W + U_d \subseteq (U+W) \cap (d+V+W).
\]
Therefore the number of representations of $d \in U-V$ as a difference of elements of $(U+W)$ and $(V+W)$ is at least $|W + U_d|$. We will need the following standard observation.

\begin{lemma}\label{PopDiff}
Let $U$ and $V$ be subsets of a commutative group and $M\geq 1$ be a parameter. For $d \in V-U$ define $V_d = V \cap (d+U)$. The set $P$ of $d \in V-U$ for which
\[
|V_d| \geq \frac{|U| |V|}{2 |U+V|}
\]
is non empty. If, moreover, 
\[
|V_d| \leq \frac{M|U| |V|}{|U+V|}
\]
for every $d \in P$, then $|P| \geq |U+V| /(2M^2)$.
\end{lemma}

\begin{proof}
By Cauchy--Schwarz,
\[
\frac{|U|^2|V|^2}{|U+V|} \leq \sum_{d \in V - U} |V_d|^2,
\]
and because
\[
\sum_{d \notin P} |V_d|^2 \leq \frac{|U||V|}{2|U+V|} \sum_d |V_d| = \frac{|U|^2 |V|^2}{2|U+V|},
\]
we get
\[
\frac{|U|^2|V|^2}{2|U+V|} \leq \sum_{d \in P} |V_d|^2.
\]
This proves $P$ is non-empty. Under the additional hypothesis we get
\[
\frac{|U|^2|V|^2}{2|U+V|} \leq \sum_{d \in P} |V_d|^2 \leq |P| \frac{M^2 |U|^2 |V|^2}{|U+V|^2}.\qedhere
\]
\end{proof}

\subsection{The Balog--Szemer\'edi--Gowers theorem}

Given two finite sets $U,V$ in a commutative group $G$ and a subset $\Gamma \subseteq U \times V$, we define the sumset of $U$ and $V$ along $\Gamma$ by
\[
U +_\Gamma V = \{ u+v : (u,v) \in \Gamma\}.
\]

The following is a version of the Balog--Szemer\'edi--Gowers theorem \cite{BaSz,Go}. It follows from a lemma in \cite{TaVu} that built on ideas in \cite{SSV}. The exact statement we use can be extracted from the proof of  \cite[Proposition 21]{MuPe}.

\begin{lemma}[Balog--Szemer\'edi--Gowers]\label{BSG}
Let $U,V,W$ be a finite sets in a commutative group $G$, let $\Gamma \subseteq U \times V$,  and $N$ be a parameter. There exist $U' \subseteq U$ and $V' \subseteq V$ of cardinalities
\[
|U'| \gg \frac{|\Gamma|}{|V|} \text{ and } |V'| \gg \frac{|\Gamma|}{|U|}
\]
such that
\[
|U'+V'| \ll \frac{|U+_\Gamma V|^3 |U|^4 |V|^4}{|\Gamma|^5}.
\]
If we further assume $|\Gamma| \geq N |V|$ and $U+_\Gamma V \subseteq W$, then there exists $\delta \gg N / |U|$ and subsets $U' \subseteq U$ and $V' \subseteq V$ of cardinalities $|U'| = \delta |U|$ and $|V'| = \delta |V|$ such that
\[
|U'+V'| \ll \frac{\delta^{-1}|W|^3 |U|^3}{N^5 |V|} |U'| \ll \frac{|W|^3 |U|^4}{N^6 |V|} |U'|.
\]
\end{lemma}

\subsection{The greedy covering lemma}

We provide a proof of a standard covering lemma because the proof of Theorem \ref{Dilates} builds on the argument presented below. 

\begin{lemma}[Greedy covering lemma]\label{GreedyCovering}
Let $U$ and $V$ be subsets of a commutative group with $|V| \geq 2$. Suppose $|U+V| \leq K |U|$. There exists a set $S \subseteq V-U$ with $|S|\ll K \log|V|$ such that $V \subseteq S+U$. There also exists a set $S' \subseteq V+U$ with $|S'|\ll K \log|V|$ such that $V \subseteq S'-U$.
\end{lemma}

\begin{proof}
Set $V_0 = V$. By Lemma \ref{PopDiff} there is some $s_1$ such that $| V \cap (s_1+U)| \geq |V|/(2K)$. Add $s_1$ to $S$ and repeat with $U$ and $V_1 := V \setminus (s_1+V)$. We still have $|U+ V_1 | \leq |U+V| \leq K |U|$ and so, by Lemma \ref{PopDiff} once again, there is some $s_2$ such that $| V_1 \cap (s_2+U)| \geq |V_1|/(2K)$. Set $V_2 = V_1 \setminus (s_2+U)$ and repeat this process until $V$ is exhausted, after, say, $\ell$ steps. 

We get $s_1, s_2, \dots, s_\ell$ and $V \supset V_1 \supset \dots \supset V_\ell = \emptyset$ with $V \setminus V_j\subseteq \{s_1 , \dots, s_j\} + U$ and
\[
|V_j| \leq \left( 1 - \frac1K\right) |V_{j-1}| \leq \dots \leq \left( 1 - \frac1K\right)^{j} |V|.
\]
After $\ell = O(K \log|V|)$ steps we get $|V_\ell| <1$ and so $V_\ell = \emptyset$. Therefore $V \subseteq \{s_1 , \dots, s_\ell\} + U$ for $\ell = O(K \log|V|)$.

A similar argument proves the existence of a suitable $S'$.
\end{proof}

\subsection{A lemma combining these tools}

A result obtained by combining the introduced tools. 

\begin{lemma}\label{Combination}
Let $U,V,W$ be a finite sets in a commutative group $G$, let $\Gamma \subseteq U \times V$,  and $N$ be a parameter. If $|\Gamma| \geq N |V|$ and $U+_\Gamma V \subseteq W$, then there exists a subset $V' \subseteq V$ of cardinality \[|V'| \gg \frac{N}{|U|} |V|\] and a set $T\subseteq G$ of cardinality
\[ |T| \ll \frac{|W|^6 |U|^8}{N^{12} |V|^2} \log|V|,\]
such that $V'+V'\subseteq T + U$.
\end{lemma}

\begin{proof}
Let $U'$ and $V'$ be the subsets given by the Balog--Szemer\'edi--Gowers result (Lemma~\ref{BSG}). Let $X$ be the subset of $U'$ that minimises the ratio $|X+V'|/|X|$ over all subsets of $U'$. The minimum ratio is at most
\[
\frac{|U'+V'|}{|U'|} \ll \frac{|W|^3 |U|^4}{N^6 |V|}.
\]
By Pl\"unnecke's inequality (Lemma \ref{Plunnecke}) we get
\[
|X + V' + V'| \ll \frac{|W|^6 |U|^8}{N^{12} |V|^2} |X|.
\]
Set $J = |W|^6 |U|^8 N^{-12} |V|^{-2}$. By the greedy covering lemma (Lemma \ref{GreedyCovering}) we can cover $V'+V'$ by $J \log|V|$ translates of $X$ and therefore by $ J \log|V|$ translates of $U$.
\end{proof}

The logarithmic factors that appear in the conclusions of Lemma \ref{GreedyCovering} and Lemma \ref{Combination}, which may be larger than $K$, will be removed at the end of the proof of Theorem~\ref{Dilates} via the tensor product trick (yet the corresponding set  does not equal any  Cartesian product).

\section{Structural results} \label{Proof 1}

In this section we derive some auxiliary results about the structure of $A$ and sets related to $A$. Broadly, we bound $A + 2\cdot A$ by partitioning $A$ into subsets $B^{(1)}, \ldots, B^{(k)}$ and applying the union bound
\[|A+2\cdot A| \leq \sum_{i=1}^k|A+2\cdot B^{(i)}|.\]
The success of the argument depends on having more information about the sumsets on the right hand side.

To that end we aim to estimate $A+2\cdot A'$ for a subset $A' \subseteq A$. Once we have done so, we will iterate the process until $A$ is exhausted. Let us begin with a basic lemma that nonetheless distinguishes between $A + 2 \cdot A'$ and $A + A' + A'$.

\begin{lemma}\label{BasicLemma}
Let $A$ be a finite subset of a commutative group $G$, let $X\subseteq A$ be such that the ratio \[K=\frac{|X+A|}{|X|}\] is minimal among all non-empty subsets of $A$, and let $A'\subseteq A$ be any subset of $A$. Suppose \[A' = \bigcup_{s \in S} A'_s\] with each $A'_s\subseteq s+X$. Then we have
 \[|A+2 \cdot A'|\leq K \sum_{s \in S} |X+A'_s|.\]
\end{lemma}

\begin{proof}
Using that $A'_s \subseteq X+s$ and Pl\"unnecke's inequality (Lemma \ref{Plunnecke}) we get
\begin{align*}
|A+2 \cdot A'| 
 & \leq \sum_{s \in S} |A + 2 \cdot A'_s| \\
 & \leq \sum_{s \in S} |A + A'_s + A'_s| \\
 & \leq \sum_{s \in S} |A + A'_s+X+s| \\ 
  & \leq K \sum_{s \in S} |X  + A'_s|. \qedhere
\end{align*}
\end{proof}

By the greedy covering lemma (Lemma~\ref{GreedyCovering}), we can choose a set $S$ of translates with $|S| \ll K \log|A|$. We also have that $|X+A_s'|\leq |X+A|$, and so
\begin{align}
\sum_{s \in S} |X+A'_s| &\leq \sum_{s\in S}|X+A| \label{sumset}\\
&\ll  K^2|X|\log|A| \label{covering}.
\end{align}
Therefore a crude upper bound for $|A+2\cdot A'|$ side is $K^3 |A|\log|A|$. However, if either of the inequalities (\ref{sumset}) or (\ref{covering}) could be improved, then we stand a chance at proving a non-trivial estimate for $|A+2\cdot A|$.

The purpose of this section is to show that $A'$ contains a large subset $B'$ for which $|A+ 2 \cdot B'|$ is smaller than $K^3 |A|$. Throughout the rest of this section, we shall use the letter $M$ for a quantity which represents some sort of saving.

Our first step is to refine the proof of the greedy covering lemma (Lemma~\ref{GreedyCovering}) and obtain more information by choosing $S$ carefully. There are similarities with \cite{Sch}. As is standard in additive combinatorics, there is a dichotomy between having more structure than expected and having a large degree of uniformity.

\begin{lemma}\label{RefinedGreedy}
Let $X$ and $A$ be subsets of an commutative group $G$ with \[|X+X|, |X+A|\leq K|X|,\] and $M \leq K$ be parameters. There exists a subset $B \subseteq A$ with  $|B| \geq |A|/3$ that satisfies at least one of the following three properties.
 \begin{itemize}
     \item[(i)] We have \[B=\bigcup_{s\in S}A_{s},\] with each 
     $A_s\subseteq s+X$, and $|S|\leq K\log|A|/M$.
     \item[(ii)] We have \[B=\bigcup_{s\in S}A_{s},\] with each 
     $A_s\subseteq s+X$, and  $|S|\leq K\log|A|$, and such that for all  $s\in S$, $|X+A_s|\leq K|X|/M$.
     \item[(iii)] If $d \in B-X$, let $B_d = B \cap (d+X)$ and let $P \subseteq B-X$ be the set
     \[P = \left\{ d \in B-X : |B_d| \geq \frac{|B||X|}{2|X+B|}\right\}.\]
    Then for each $d\in P$,
     \begin{enumerate}
         \item $|B_d|\leq M|B||X|/|X+B|$, and
         \item $\max\{|X+X|, |X+B|\} / M \leq |X+B_d|$.
     \end{enumerate} 
 \end{itemize}
\end{lemma}

In this lemma, cases (\emph{i}) and (\emph{ii}) can be interpreted as providing additional structure, and will be useful in improving the inequalities (\ref{sumset}) and (\ref{covering}), respectively. Case (\emph{iii}) grants us uniformity, which will be exploited in Lemma \ref{TechnicalLemma}.
\begin{proof}
We will cover $A$ by translates of $X$, as in the proof of the greedy covering lemma (Lemma~\ref{GreedyCovering}). In fact, we will show $A\subseteq (S^{(1)} \cup S^{(2)} \cup S^{(3)}) + X$ for some sets $S^{(1)}$, $S^{(2)}$, and $S^{(3)}$ that will be defined in the process. There is asymmetry in the role of $S^{(i)}$. $B$ is determined by $A\cap(S^{(1)}+X)$ or $A\cap(S^{(2)}+X)$ (if either of them is large) or by $A \cap (s+ X)$ for a single element of $S^{(3)}$ (if both $A\cap(S^{(1)}+X)$ or $A\cap(S^{(2)}+X)$ are small).

First, we construct a decreasing sequence of sets $B^{(j)}$, a sequence of elements $s_j\in B^{(j)}$, and a sequence of sets $A_{s_j}\subseteq B^{(j)}$ for $j\geq 0$. To initialize set $B^{(0)}=A$. Next, suppose $B^{(j)} \subseteq A$ has been defined. For $d \in B^{(j)} - X$ set \[B^{(j)}_d = B^{(j)} \cap (d +X).\] Let \[P^{(j)}=\left\{d:|B^{(j)}_d| \geq \frac{|B^{(j)}||X|}{2|X+B^{(j)}|}\right\}.\] 
That $P^{(j)}$ is non-empty follows from Lemma~\ref{PopDiff}. 

Suppose $d_*\in P^{(j)}$ maximizes $|B^{(j)}_d|$ and suppose $d_{**}\in P^{(j)}$ minimizes the quantity $|X+B^{(j)}_d|$.

To define $s_j,\ A_{s_j}$ and $B^{(j+1)}$ we examine three cases.

\emph{Case 1}: If \[|B^{(j)}_{d_*}| \geq \frac{M |B^{(j)}||X|}{|X+B^{(j)}|}\] we then  define  
\[s_j=d_*,\ A_{s_j}=B^{(j)}_{d_*},\ B^{(j+1)} = B^{(j)} \setminus B^{(j)}_{d_*},\]
and add $d_*$ to $S^{(1)}$. We also note that since $B^{(j)}_{d_*}\subseteq B^{(j)}\subseteq A$, we have
\begin{equation} \label{S_1}
    \frac{|B^{(j+1)}|}{|B^{(j)}|} \leq 1-\frac{Mu}{K}.
\end{equation}
This is because 
\begin{align*}
    \frac{|B^{(j+1)}|}{|B^{(j)}|} 
    & = 1-\frac{|B^{(j)}_{d_*}|}{|B^{(j)}|} \\
    & \leq 1-\frac{M |B^{(j)}||X|}{|X+B^{(j)}||B^{(j)}|} \\
    & \leq 1-\frac{M|X|}{|X+A|} \\
    & \leq 1-\frac{M}{K}.
\end{align*}
\emph{Case 2}: If 
\[|X+B^{(j)}_{d_{**}}| \leq \frac{\max\{|X+X|, |X+B|\}}{M}\] we then define  
\[s_j=d_{**},\ A_{s_j}=B^{(j)}_{d_{**}},\ B^{(j+1)} = B^{(j)} \setminus B^{(j)}_{d_{**}},\]
and add $d_{**}$ to $S^{(2)}$. We also note that since $B^{(j)}_{d_{**}}\subseteq B^{(j)}\subseteq A$ and ${d_{**}}\in P^{(j)}$,
\begin{equation} \label{S_2}
    \frac{|B^{(j+1)}|}{|B^{(j)}|}\leq \left(1 -\frac{|B^{(j)}_{d_{**}}|}{|B^{(j)}|} \right) \leq \left( 1 -\frac{1}{2K} \right)
\end{equation}
and
\begin{equation} \label{S_2'}
    |X+A_{s_j}| \leq \frac{K}{M} |X|.
\end{equation}

\emph{Case 3}: If cases 1 and 2 fail to apply, then for every $d \in P^{(j)}$ we have \[|X+B^{(j)}_{d}| \geq \frac{\max\{ |X+X|, |X+B|\}}{M}.\] As in case 2, we define
\[s_j=d_{**},\ A_{s_j}=B^{(j)}_{d_{**}},\ B^{(j+1)} = B^{(j)} \setminus B^{(j)}_{d_{**}},\] but add $d_{**}$ to $S^{(3)}$. This time note that $B^{(j)}$ satisfies (\emph{iii}). 

We repeat the process for $B^{(j+1)}$ until we exhaust $A$ and get $B^{(\ell+1)} = \emptyset$. We obtain $S= S^{(1)} \cup S^{(2)} \cup S^{(3)}$ and
\[
A = \bigcup_{s \in S} A_s.
\]
Note here that 
\[
\frac{1}{|A|} \leq \frac{|B^{(\ell)}|}{|A|} = \frac{|B^{(\ell)}|}{|B^{({\ell-1})}|} \dots \frac{|B^{(1)}|}{|B^{(0)}|}.
\]
It follows from the above and \eqref{S_1} and \eqref{S_2} that $|S^{(1)}| \leq K \log|A| / M$ and $|S^{(2)}| \leq K \log|A|$.

For $i=1,2$ define
\[
A^{(i)} = \bigcup_{s \in S^{(i)}} A_{s}.
\]

If $|A^{(1)}| \geq |A|/3$ we set $B = A^{(1)}$ and observe that condition (\emph{i}) is satisfied, while if $|A^{(2)}| \geq |A|/3$ we set $B = A^{(2)}$ and observe that condition (\emph{ii}) is satisfied. Otherwise, let $j_*$ be the first $j$ with $s_j\in S^{(3)}$ and set $B = B^{(j_*)}$. We have seen that $B$ satisfies the conditions in (\emph{iii}) and moreover, \[|B| \geq |A| - |A^{(1)}| - |A^{(2)}| \geq |A|/3.\qedhere\] 
\end{proof}

Our next task is to turn the uniformity case, (\emph{iii}) of Lemma \ref{RefinedGreedy}, into structural information. We will use ideas from \cite{KaKo} and \cite{ScSh} (see also \cite{Sch,Shk1,Shk2,Shk3}). In a nutshell, we will use uniformity to obtain an estimate for the additive energy of the set $P$ which depends on $M$, but crucially does not depend on $K$. This means that $P$ has considerably more additive structure than $A$, and subsequent methods from additive combinatorics are much less costly. It is crucial for our method that the sets $X+X$, $X+B$, $P$, and all the $X+B_d$ have, up to powers of $M$, the same cardinality.

\begin{lemma}\label{TechnicalLemma}
Let $A$ and $B$ and $X$ be finite subsets of a commutative group $G$ and $M \leq K \leq |X|$ be parameters such that $|X+X| , |X+B| \leq K |X|$. If $B$ and $X$ satisfy (iii) of Lemma \ref{RefinedGreedy}, then there exists a subset $B ' \subseteq B$ of cardinality 
\[
|B'| \gg \frac{|B|}{M^{3}}
\]
such that 
\[
|A + B' + B'| \ll  M^{16} |X+X+A|\log|X|.
\]
\end{lemma}

\begin{proof}
Let us denote 
\[N = \min\{|X+B_d| : d \in P\}.\]
Then (\emph{b}) of Lemma \ref{RefinedGreedy} tells us that
\begin{equation}\label{N bound}\frac{\max\{|X+B|,|X+X|\}}{M}\leq N.\end{equation}
In addition the observation from \cite{KaKo} that for all $d \in B-X$
\[X + B_d \subseteq (X+B) \cap (d+X+X),\]
from which we deduce the estimate
\[N\leq |(X+B) \cap (d+X+X)|\]
for all $d\in P$. Taking the sum of this inequality over all $d\in P$ yields
\begin{equation}\label{Popular}
N|P|\leq \sum_{d \in P} |(X + B)\cap(d+X+X)| = \sum_{z \in X+B} |(X+X) \cap (z-P)|.
\end{equation}
We define $\Gamma \subseteq (X+X) \times P$ by
\[
\Gamma = \{ (u, v) \in (X+X) \times P : u+v \in X+B\},
\]
so that \[(X+X)+_\Gamma P \subseteq X+B.\] The inequality (\ref{Popular}) shows that $|\Gamma| \geq N |P|$ by simple double counting. 

Next, by Lemma \ref{PopDiff} and by (\emph{a}) of Lemma \ref{RefinedGreedy},  we get
\begin{equation}\label{P bound}
\frac{|X+B|}{2M^2}\leq |P|.
\end{equation}

We apply Lemma~\ref{Combination} to $U = X+X$, $V = P$ and $W = X+B$ yielding sets $P'\subseteq P$ and $T$. The cardinality of $P'$ is guaranteed to satisfy
\[|P'| \gg \frac{N}{|X+X|} |P| \geq \frac{|P|}{M},\]
by (\ref{N bound}). Moreover, $P'+P' \subseteq T + (X+X)$ and
\[
|T| \ll \frac{|X+B|^6 |X+X|^8}{N^{12} |P|^2} \log|P| \ll M^{16} \log|P|,
\]
the final estimate being a consequence of (\ref{P bound}) and (\ref{N bound}).

We change task now and locate a translate of a subset of $B$ inside $P'$. From the definition of $P$ and the fact that $P'\subseteq P$,
\[
\sum_{x \in X} |(B - x) \cap P'| \geq \frac{|P'| |B| |X|}{2|B+X|} \gg \frac{|P| |B||X|}{M |B+X|}\gg \frac{|X| |B|}{M^{3}}, 
\]
where we have used (\ref{P bound}) in the final inequality.
From this, there exists an $x_0 \in X$ such that 
\[
|B \cap (x_0+P')| \gg \frac{|B|}{M^3}.
\] 
Set \[B' = B \cap (x_0 + P').\]
We have
\begin{align*}
    B'+B'+A
    \subseteq  2x_0+P'+P' +A 
     \subseteq 2x_0+T+ X+X +A,
\end{align*}
from which we get the desired \[|A+ B' + B'|\leq |T||X+X+A|\ll  M^{16} |X+X+A| \log|P|.\qedhere\]
\end{proof}

We now have some refined information in both the structured and uniform cases. The main lemma of this section combines Lemma \ref{RefinedGreedy} and Lemma \ref{TechnicalLemma} to produce the dichotomy behind proof of Theorem~\ref{Dilates}: given a subset $A'$ of $A$, either $A'$ has a large subset $B$ such that $|A + 2 \cdot B|$ is a little smaller than $K^3 |A|$, or else $A'$ has a slightly smaller subset $B$ such that $|A + 2 \cdot B|$ is only a little greater than $K^2 |A|$.

\begin{lemma}[Main lemma]\label{MainLemma}
Let $A$ be a finite subset of a commutative group $G$ with the property that $|A+A| \leq K |A|$, let $X\subseteq A$ be the subset that minimises the ratio $|X+A|/|X|$, and let $A'\subseteq A$ be any subset of $A$. For every $M \leq K$, there exists a subset $B \subseteq A'$ which satisfies one of the following two pairs of properties. 
\begin{enumerate}
\item $|B| \gg |A'|$ and \[|A+ 2 \cdot B| \leq \frac{K^3|X|\log|A|}{M}.\]
     \item $|B| \gg |A'|/M^{3}$ and \[|A+ 2 \cdot B| \ll M^{16}K^2|A|\log|A| .\]
 \end{enumerate}
\end{lemma}

\begin{proof}
We apply Lemma \ref{RefinedGreedy} to find a subset $B_*$ with $|B_*|\gg |A'|$. If $B_*$ satisfies (\emph i) of Lemma \ref{RefinedGreedy}, then set $B = B_*$ and note that by Lemma  \ref{BasicLemma}, 
\[
|A+ 2\cdot B| \leq K \sum_{s \in S} |A'_s + X| \leq K\cdot|S|\cdot K |X| \leq \frac{K^3 \log|A|}{M} |X|. 
\]
If $B_*$ satisfies (\emph{ii}) of Lemma \ref{RefinedGreedy}, then set $B = B_*$ and note that by Lemma ~\ref{BasicLemma}, 
\[
|A+ 2\cdot B| \leq K \sum_{s \in S} |A'_s + X| \leq K \cdot |S| \cdot \frac{K}{M} |X| \leq \frac{K^3 \log|A|}{M} |X|. 
\]

If $B_*$ instead satisfies (\emph{iii}), then we apply Lemma \ref{TechnicalLemma} to $B_*$ to find a subset $B$ such that \[|B|\gg \frac{|A'|}{M^{3}}\] and note that 
\[
|A+2 \cdot B| \leq |A+B+B| \ll M^{16}K^2 |A|\log|A|,
\]
because by Pl\"unnecke's inequality (Lemma~\ref{Plunnecke}) we have 
\[
|X+X+A| \leq |X+A+A| \leq K^2 |X| \leq K^2 |A|. \qedhere
\]
\end{proof}

\section[Proof of Theorem 1]{Proof of Theorem~\ref{Dilates}} \label{Proof2}

In this section we prove Theorem \ref{Dilates} by iterated application of Lemma \ref{MainLemma}. We partition $A$ into subsets $B^{(1)}, \ldots, B^{(k)}$ as follows. First, we apply Lemma~\ref{MainLemma} with $A'=A$ to find a subset $B^{(1)}$. We then set $A' = A \setminus B^{(1)}$ are repeat the process until $A$ is exhausted. From the union bound, we have
\begin{equation}\label{UnionBound}|A+2\cdot A| \leq \sum_{i=1}^k|A+2\cdot B^{(i)}|.\end{equation}
We are left to estimate the right hand side.

In the iteration, the number of sets $B^{(j)}$ coming from alternative (\emph{a}) of Lemma \ref{MainLemma} can be at most $O(\log|A|)$. Indeed, each time this case occurs, the set $A'\setminus B^{(j)}$ is smaller by a constant factor. Thus such sets contribute at most $O(K^3 (\log|A|)^2 |A| / M)$ to the right hand side of (\ref{UnionBound}). 

In a similar way, alternative (\emph{b}) will exhaust $A$ after occurring $O(M^3 \log|A|)$ times, and therefore the contribution to the right hand side of (\ref{UnionBound}) that comes from such $B^{(j)}$ is $O(K^2 M^{19} (\log|A|)^2 |A|)$. 

From our analysis, we conclude
\[
|A+ 2 \cdot A| \ll \left( \frac{K^3 (\log|A|)^2}{M} + K^2 M^{19} (\log|A|)^2  \right) |A|.
\]
Choosing $M = K^{1/20}$ gives
\begin{equation} \label{LogDilates}
|A+ 2 \cdot A| \leq C (\log|A|)^2 K^{3- 1/20} |A|,
\end{equation}
for some absolute constant $C>0$.

Our final task is to remove $C$ and the logarithmic terms. We have proved an upper bound on $A + 2 \cdot A$ that holds for any finite set $A$ of an arbitrary commutative group $G$. Given such a pair $(A, G)$ we apply \eqref{LogDilates} to the $r$-fold product $A^r$ in $G^r$. Note that $|A^r| = |A|^r$, $|A^r + A^r| = |A+A|^r$ (and therefore $K$ becomes $K^r$) and $|A^r + 2 \cdot A^r| = |A+2 \cdot A|^r$. Applying \eqref{LogDilates} gives
\[
|A + 2 \cdot A|^r \leq C r \log|A| (K^{3 - 1/20})^r |A|^r.
\]
Taking $r$-th roots gives
\[
|A + 2 \cdot A| \leq C^{1/r} r^{1/r} \log|A|^{1/r} K^{3 - 1/20} |A|.
\]
Letting $r \to \infty$ finishes the proof of Theorem~\ref{Dilates}.

\section{An upper bound on $|A-2 \cdot A|$} \label{A-2A}

Theorem~\ref{Dilates} can be generalised to the set
\[
A - 2 \cdot A = \{a - a'-a' : a,a' \in A\}
\]
under the additional assumption that $|A-A|$ is bounded above by $K|A|$. This is not a particularly restrictive condition. For example in \cite{Bu1,Bu2} the bound used is on $|A-A|$, but it comes from an application of the Balog--Szemer\'edi--Gowers theorem in \cite{SSV} and so an identical bound can be proved for $|A+A|$. See \cite[Lemma 6]{Bu1} for details. Similarly, in \cite{CrSi} the bound on $|A+A|$ comes from $A$ being an interval in $\mathbb{Z}$. Therefore an identical upper bound can be proved for $|A-A|$ (the same is true for subsets of arithmetic progressions of large relative density). See \cite[Corollary 5.4]{CrSi} for details.

\begin{theorem}\label{Dilates diff}
Let $A$ be a subset of a commutative group $G$ and $K$ be a parameter. If $\max\{|A+A|, |A-A|\} \leq K |A|$, then $|A - 2 \cdot A| \leq K^{2.95} |A|.$
\end{theorem}

\begin{proof}[Sketch of proof]
We carry out a similar argument to that of the proof of Theorem~\ref{Dilates}.

In Lemma~\ref{BasicLemma}, we cover $A'$ by translates of $-X$ using Lemma \ref{GreedyCovering}: $A' \subseteq S' - X$. This gives
\[
A' = \bigcup_{s \in S'} A'_s \text{ where } A'_s = A' \cap (s-X).
\]
Repeating the steps in the proof of Lemma ~\ref{BasicLemma} gives
\[
|A - 2\cdot A'| \leq \sum_{s \in S'} |A-A'_s-(s-X)|\leq K \sum_{s \in S'} |A'_s-X|.
\]

In Lemma \ref{RefinedGreedy} we assume $|X-A| \leq K |A|$. We then take $s\in B+X$ and $A_s \subseteq B \cap ( s - X)$. In part (\emph{ii}) we require $|A_s - X| \leq K|X| / M$. In part (\emph{iii}) we consider popular sums in $X+B$ (the popularity parameter remains $|X||B| / (2 |X+B|)$) and replace $X+B$ by $B-X$.

In Lemma \ref{TechnicalLemma} we make the natural adjustment when applying the Balog--Szemer\'edi--Gowers result: $U = -(X+X), V=P, W=B-X$. This gives a subset $U'$ and $P'$ such that $|U' + P'| \leq M^8 |U'|$ and eventually a subset $B'$ that is covered by $M^{16} \log|P|$ translates of $-(X+X)$. Hence $A - B' - B'$ is covered by $M^{16} \log|P|$ translates of $A+X+X$. The remainder of the proof remains largely unchanged.
\end{proof}

Our method does not work under the hypothesis $|A+A| \leq K |A|$ on its own. When we apply the Balog--Szemer\'edi--Gowers theorem, we need three sets that have nearly equal cardinality.

Returning to the hypothesis $|A+A| \leq K |A|$, the authors are not aware of examples where $|A+ 2 \cdot A| \geq K^2 |A|$; there are, however, arbitrarily large sets $A$ such that $|A- 2 \cdot A|$ is  far larger than $K^{5/2} |A|$. The examples are essentially due to Fre\u{\i}man and Pigarev \cite{FrPi} (analysed in \cite{Gr,HRY}).  
\begin{example}
Let 
\[
q=  \frac{2\log(1 + \sqrt2)}{\log2} > 2.543.
\]
There exists an infinite family of sets $A$ in commutative groups such that if we set $K = |A+A| / |A|$, then
\[
|A- 2 \cdot A| \geq \frac{K^{q}}{2 \log_2K} |A|.
\]
\end{example}

\begin{proof}
We think of $d$ as being a fixed positive integer and $T$ an integer tending to infinity. Set
\[
A = \{x \in \mathbb{Z}^d : x_i \geq 0 \text{ for all } i \,,\sum_{i=1}^d x_i \leq T \}.
\]

A basic counting argument shows that $|A| = \binom{T+d}{d}$ is asymptotically equal to $T^d / d!$. To find $K$ note that
\[
A+A = \{x \in \mathbb{Z}^d : x_i \geq 0 \text{ for all } i \,, \sum_{i=1}^d x_i \leq 2T \}.
\]
So $|A+A| = \binom{2T+d}{d}$, which makes $K$ asymptotically equal to $2^d$. Now
\[
A-2 \cdot A \subseteq \{x \in \mathbb{Z}^d : -2 T \leq \sum_{i=1}^d x_i \leq T \}.
\]
To estimate the cardinality of $A- 2 \cdot A$ we denote by $p,z,n$ the number positive, zero, and negative coordinates. There are $\binom{d}{p,z,n}$ partitions of indices for a given triple $(p,z,n)$. Given a triple, the sum of the negative coordinates is at least $-2T$ and so there are $\binom{2T}{n}$ possibilities. The sum of the  positive coordinates is at most $T$ and so there are $\binom{T}{p}$ possibilities. Therefore
\begin{align*}
    |A- 2 \cdot A| 
    & = \sum_{p+z+n=d} \binom{d}{p,z,n} \binom{2T}{n} \binom{T}{p}\\
    & \geq \sum_{n=0}^d \binom{d}{n} \binom{2T}{n} \binom{T}{d-n} \\
    & = (1+o(1)) \frac{T^d}{d!} \sum_{n=0}^d \binom{d}{n}^2 2^n\\ 
    & = (1+o(1)) |A| \sum_{n=0}^d \binom{d}{n}^2 2^n \\
    & \geq (1+o(1)) |A| \frac{\left( \sum_{n=0}^d \binom{d}{n} 2^{n/2}\right)^2}{d+1} \\
    & = (1+o(1)) |A| \frac{(1+\sqrt{2})^{2d}}{d+1}.
\end{align*}

By the definition of $q$, $2^{q/2} = (1 + \sqrt2)$, and $d=(1+o(1))\log_2K$, so that 
\[
\frac{(1+\sqrt2)^{2d}}{d+1} \geq \frac{2^{q d}}{d+1} = (1+o(1)) \frac{K^{q}}{1+\log_2K}.
\]
Therefore
\[
|A - 2 \cdot A| \geq (1+o(1)) \frac{K^{q}}{1+\log_2K} |A|. \qedhere
\]
\end{proof}

\section{Upper bounds on $|A+\lambda \cdot A|$ for small prime $\lambda$} \label{A+ prime A}

In this section we denote by $\lambda * A$ the $\lambda$-fold sumset of $A$:
\[
\lambda * A = \{ a_1 + \dots + a_\lambda : a_1, \dots, a_\lambda \in A\}.
\]
Dilate sums for integer $\lambda$ are defined in the natural way. The proof of Theorem~\ref{Dilates} works for all $\lambda >0$ and yields the following result.

\begin{theorem} \label{Large dilates}
Let $\lambda>0$ be a positive integer. Set
\begin{equation} \label{c lambda}
c_\lambda = \frac{\lambda-1}{4 + 8 \lambda} = \frac{1}{8} - \frac{3}{8(1+2 \lambda)}.
\end{equation}
For every subset $A$ of a commutative group $G$ and every $K \in \mathbb{R}$ such that $|A+A| \leq K |A|$, we have $|A+ \lambda \cdot A| \leq K^{\lambda+1 - c_\lambda} |A|.$
\end{theorem}

\begin{proof}[Sketch of proof]
In Lemma~\ref{BasicLemma} we cover $A + \lambda \cdot A'$ by the union of the $A+ \lambda \cdot A'_s$. For each $s$, we bound $|A+ \lambda \cdot A'_s|$ by a combination of inclusions and Lemma~\ref{Plunnecke}:
\[
|A+ \lambda \cdot A'_s| \leq |X+ A'_s +(\lambda-2)*A| \leq K^{\lambda-2} |X+A'_s|.
\]
Lemma~\ref{RefinedGreedy} remains unchanged. In Lemma~\ref{TechnicalLemma} we cover $\lambda * B'$ by $M^{8 \lambda} \log|P|$ translates of $X+X$. As a result in Lemma~\ref{MainLemma} we balance, up to logarithms, the terms
\[
\frac{K^{\lambda+1}}{M} |A| \quad \text{ and } \quad M^{3 + 8 \lambda} K^2 |A|.
\]
We set $M = K^{c_\lambda}$ and apply the tensor power trick.
\end{proof}

The bound improves on the Pl\"unnecke-type bound that comes from $A + \lambda A \subseteq (\lambda+1) * A$. We investigate for which $\lambda$ Theorem~\ref{Large dilates} may represent the state-of-the-art. We begin by refining some observations of Bukh \cite{Bu2}.

\begin{lemma}\label{dilate sums}
Let $A$ a be a set in a commutative group, $K = |A+A| / |A|$  and $X \subseteq A$ be the subset that minimises $|X+A| / |X|$.
\begin{enumerate}
    \item For all integers $1 \leq \lambda_1, \lambda_2$ \[ |A+ (\lambda_1 \pm \lambda_2) \cdot A| \le K \frac{|X+ \lambda_1 \cdot A| \, |X+ \lambda_2 \cdot A|}{|X|}\le K^{\lambda_1+\lambda_2 +1} |A|.\]
    \item For all integer $\lambda_1, \lambda_2 \geq 1$ \[ |A \pm (\lambda_1 \lambda_2) \cdot A| \le \frac{|A + \lambda_1 \cdot X| \, | X + \lambda_2 \cdot A|}{|X|} \le K^{\lambda_1+\lambda_2} |A|. \] 
    \item For all integer $\lambda \geq 1$ and $j \ge 2$ \[ |A \pm \lambda^{j} \cdot A| \le \left( \frac{|A + \lambda \cdot X|}{|X|} \right) \left( \frac{|X + \lambda \cdot X|}{|X|} \right)^{j-2} |X + \lambda \cdot A| \le K^{j\lambda} |A|.\] 
\end{enumerate}
\end{lemma}

\begin{proof}
We make repeated use of Pl\"unnecke's inequality (Lemma \ref{Plunnecke}) and of the following combined corollary of Pl\"unnecke's inequality and Ruzsa's triangle inequality \cite{Ru1}: For all sets $U, V, W$ we have
\[
|U\pm V| \le \frac{|U+W| \, |V+W|}{|W|}.
\]

The first inequality in the first claim follows immediately from this and Lemma~\ref{Plunnecke}:
\begin{align*}
        |A+ (\lambda_1 \pm \lambda_2) \cdot A| \le\frac{|X+ A + \lambda_1 \cdot A| \, |X+ \lambda_2 \cdot A|}{|X|} \leq K \frac{|X+ \lambda_1 \cdot A| \, |X+ \lambda_2 \cdot A|}{|X|}.
\end{align*}
Lemma~\ref{Plunnecke} once more gives $|X + \lambda \cdot A| \leq |X + \lambda * A| \leq K^\lambda |X|.$

For the second claim note that
\[
|A \pm (\lambda_1 \lambda_2) \cdot A| \le \frac{|A + \lambda_1 \cdot X| \, |\lambda_1 \cdot X + (\lambda_1 \lambda_2) \cdot A|}{|\lambda_1 \cdot X|} = \frac{|A + \lambda_1 \cdot X| \, | X + \lambda_2 \cdot A|}{|X|}.
\]
Inclusion in sumsets and Lemma~\ref{Plunnecke} give the second inequality.

To prove the third claim, observe that for all $\lambda, j \geq 1$, by Lemma~\ref{Plunnecke}
\[
|X + \lambda^j \cdot A| \leq \frac{|X+\lambda \cdot X| \, |\lambda \cdot X + \lambda^j \cdot A|}{|\lambda \cdot X|} = \frac{|X+\lambda \cdot X|\, | X + \lambda^{j-1} \cdot A|}{|X|}. 
\]
By induction for all $\lambda, j \geq 1$
\[
|X + \lambda^j \cdot A| \leq \left(\frac{|X+\lambda \cdot X|}{|X|} \right)^{j-1} | X + \lambda \cdot A|. 
\]
To prove the first inequality in the third claim for $\lambda \ge 1$ and  $j \geq 2$ we apply Lemma~\ref{Plunnecke} and the above
\begin{align*}
|A \pm \lambda^{j} \cdot A| 
& \leq \frac{|A + \lambda \cdot X|\, |\lambda \cdot X + \lambda^{j} \cdot A|}{|\lambda \cdot X|} \\
& = \frac{|A + \lambda \cdot X| \, |X + \lambda^{j-1} \cdot A|}{|X|}\\
& \le \left( \frac{|A + \lambda \cdot X|}{|X|} \right) \left( \frac{|X + \lambda \cdot X|}{|X|} \right)^{j-2} |X + \lambda \cdot A|.
\end{align*}
We apply Lemma~\ref{Plunnecke} again to bound this by $K^{j \lambda} |A|$.
\end{proof}

We list what appear to be the known upper bounds on $|A + \lambda \cdot A|$ for small positive integers $\lambda$. We denote by $p_\lambda$ the infimum of permissible exponents:
\[
p_\lambda = \inf\{ p : |A+ \lambda \cdot A| \leq K^p |A| \text{ for all $A$}\};
\]
and recall the definition of the  $c_\lambda$ in \eqref{c lambda}.
\begin{itemize}
    \item $\lambda=2$: $p_2 \leq 3 - c_2 = 3-1/20$ (Theorem~\ref{Large dilates}).
    \item $\lambda=3$: $p_3 \leq 4 -c_3 = 4 - 1/14$ (Theorem~\ref{Large dilates}).
    \item $\lambda=4$: $p_4 \leq 4$ ($\lambda=j=2$ in the third part of Lemma~\ref{dilate sums}).
    \item $\lambda = 5$: $p_5  \leq  6 - c_5  = 6 - 1/11$ (Theorem~\ref{Large dilates}).
    \item $\lambda=6$: $p_6 \leq 5$ ($\lambda_1=2$ and $\lambda_2=3$ in the second part of Lemma~\ref{dilate sums}).
    \item $\lambda=7$: $p_7 \leq 7$ ($\lambda_1=1$ and $\lambda_2=6$ in the first part of Lemma~\ref{dilate sums}).
\end{itemize}

We stop here because for $\lambda > 5$ it seems that $p_\lambda$ is never close to $\lambda+1$ and so using Theorem \ref{Large dilates} is not optimal.

\section{A sharper upper bound for large $K$} \label{Final Section}

When $K$ is a large power of $|A|$ Pl\"unnecke's inequality seizes to be the best known bound on $|A+A+A|$. Ruzsa investigated this question in~\cite[Section 6]{Ru3} (see also~\cite[Section 1.9]{Ru4}). For example,~\cite[Theorem 9.1]{Ru4} implies
\[
|A+A+A| \leq K^2 |A|^{3/2}.
\]
This is stronger than the $K^3|A|$ bound that comes from Pl\"unnecke's inequality when $K\geq |A|^{1/2}$. For $|A+2 \cdot A|$, however, the induced bound in this range is worse than the trivial bound $|A+2 \cdot A| \leq |A|^2$. In this final section we prove \eqref{Large K}.
\begin{theorem}\label{Dilates large K}
Let $A$ be a subset of a commutative group $G$ and $K$ be a parameter. If $|A+A| \leq K |A|$, then $|A+2 \cdot A| \leq (K |A|)^{4/3}.$
\end{theorem}

Theorem~\ref{Dilates large K} is non-trivial and improves Theorem~\ref{Dilates} when $|A|^{20/97} \leq K \leq |A|^{1/2}$. The proof, a variant of arguments of Ruzsa~\cite{Ru3,Ru4}, is based on Pl\"unnecke's inequality for a large subset \cite{KaSh,Ru4}, the trivial bound $|Z + 2 \cdot Z| \leq |Z|^2$, and the tensor power trick. Let us first state~\cite[Corollary 1.7.5]{Ru4}.  

\begin{lemma}[Pl\"unnecke's inequality for a large subset] \label{Large Plu}
Let $A$ be a finite subset of a commutative group, and $K>0$ and $0 < \delta < 1$ be positive real numbers. If $|A+A| \leq K|A|$, then there exists a subset $Y \subseteq A$ such that $|Y| \geq (1-\delta) |A|$ and 
\[
|Y+A+A| \leq \frac{1}{\delta^2} K^2 |Y| - \frac{(1-\delta) (2-\delta)}{2 \delta^2} K^2 |A|.
\] 
\end{lemma}

\begin{proof}[Proof of Theorem~\ref{Dilates large K}]
We may suppose that $K < |A|^{1/2}$ (otherwise  $|A+ 2 \cdot A| \leq |A|^2 \leq (K |A|)^{4/3}$). We apply Lemma~\ref{Large Plu} for $\delta = K^{2/3} / |A|^{1/3} < 1$. The reasons behind this choice of $\delta$ become apparent further down. We set $Z = A \setminus Y$ to be the complement of $Y$ in $A$ and observe that
\[
A+ 2 \cdot A \subseteq (Y+A+A) \cup (Z+2\cdot Z).
\]
Therefore
\begin{equation} \label{LK1}
|A+2 \cdot A| \leq |Y+A+A| + |Z+2\cdot Z| \leq \frac{1}{\delta^2} K^2 |Y| - \frac{(1-\delta) (2-\delta)}{2 \delta^2} K^2 |A| + (|A|-|Y|)^2.
\end{equation}

The right side of~\eqref{LK1} is a monic quadratic polynomial in $|Y|$ and therefore attains its maximum value at the endpoints of the domain. When $|Y| = (1-\delta)|A|$, \eqref{LK1} becomes
\begin{equation} \label{LK2}
|A+2 \cdot A| \leq \frac{1-\delta}{2\delta} K^2 |A| + \delta^2 |A|^2 \leq \frac{1}{2\delta} K^2 |A| + \delta^2 |A|^2.
\end{equation}
While, when $|Y| = |A|$, \eqref{LK1} becomes
\begin{equation} \label{LK3}
|A+2 \cdot A| \leq \frac{3-\delta}{2\delta} K^2 |A| \leq \frac{3}{2\delta} K^2 |A|.
\end{equation}

Our choice of $\delta = K^{2/3} / |A|^{1/3}$ makes the right side of both ~\eqref{LK2} and~\eqref{LK3} equal to $3 (K|A|)^{4/3}/2$. An application of the tensor product trick removes the $3/2$ and proves the theorem.
\end{proof}

For every positive integer $\lambda$, a variation of the above gives
\[
|A + \lambda \cdot A| \leq (K |A|)^{\tfrac{2\lambda}{\lambda+1}}.
\] 
This is non-trivial and improves both Bukh's bound~\cite{Bu2} and those described in Section~\ref{A+ prime A} for small enough $\lambda$ and large enough $K$ (in terms of $\lambda$ and $|A|$).


\section*{Acknowledgments} 
The content and exposition in this article has benefited from insightful discussions with Boris Bukh, Brendan Murphy, Oriol Serra, and Ilya Shkredov. The authors also thank the anonymous referee for a careful reading of the paper and helpful suggestions. 

\bibliographystyle{amsplain}

\begin{thebibliography}{99}

\bibitem[BaSh]{BaSh}
A. Balog and G. Shakan.
\newblock On the sum of dilations of a set.
\newblock {\em Acta Arith.}, 164(2):153--162, 2014. 

\bibitem[BaSz]{BaSz}
A. Balog and E. Szemer\'edi. 
\newblock A statistical theorem of set addition.
\newblock {\em  Combinatorica}, 14(3):263--268, 1994. 

\bibitem[BrGr]{BrGr}
E. Breuillard and B. Green. 
\newblock Contractions and expansion.
\newblock {\em  European J. Combin.}, 34(8):1293--1296, 2013. 

\bibitem[Bu1]{Bu1}
B. Bukh. 
\newblock Non-trivial solutions to a linear equation in integers.
\newblock {\em Acta Arith.}, 131(1):51--55, 2008. 

\bibitem[Bu2]{Bu2}
B. Bukh. 
\newblock Sums of dilates.
\newblock {\em  Combin. Probab. Comput.}, 17(5):627--639, 2008.

\bibitem[BuZh]{BuZh}
A. Bush and Y. Zhao. 
\newblock New upper bound for sums of dilates.
\newblock {\em  Electron. J. Combin.}, 24(3):7pp, 2017.

\bibitem[ChFa]{ChFa}
Y.-G. Chen and J.-H. Fang. 
\newblock Sums of dilates in the real numbers.
\newblock {\em Acta Arith.}, 182(3):231--241, 2018. 

\bibitem[CHS]{CHS}
J. Cilleruelo, Y. O. Hamidoune and O. Serra. 
\newblock On sums of dilates.
\newblock {\em  Combin. Probab. Comput.}, 18(6):871--880, 2009.

\bibitem[CSV]{CSV}
J. Cilleruelo, M. Silva, and C. Vinuesa. 
\newblock A sumset problem.
\newblock {\em  J. Comb. Number Theory}, 2(1):79--89, 2010. 

\bibitem[CrSi]{CrSi}
E. Croot and O. Sisask. 
\newblock A probabilistic technique for finding almost-periods of convolutions.
\newblock {\em  Geom. Funct. Anal.}, 20(6):1367--1396, 2010.

\bibitem[DCS]{DCS}
S.-S. Du, H.-Q. Cao, Z.-W. Sun,
\newblock On a sumset problem for integers.
\newblock {\em  Electron. J. Combin.}, 21(1):25pp, 2014.

\bibitem[F-P]{F-P}
G. Fiz Pontiveros. 
\newblock Sums of dilates in $\mathbb{Z}_p$.
\newblock {\em  Combin. Probab. Comput.}, 22(2):282--293, 2013.

\bibitem[Fr]{Fr}
G. A. Fre\u{\i}man. 
\newblock Foundations of a structural theory of set addition.
\newblock Translations of Mathematical Monographs, Vol 37. American Mathematical Society, Providence, R. I., 1973. vii+108 pp.. 

\bibitem[FrPi]{FrPi}
G. A. Fre\u{\i}man and V. P. Pigarev. 
\newblock The relation between the invariants $R$ and $T$.
\newblock Number-theoretic studies in the Markov spectrum and in the structural theory of set addition (Russian), pp. 172--174. Kalinin. Gos. Univ., Moscow, 1973. 

\bibitem[Go]{Go}
W. T. Gowers. 
\newblock A new proof of Szemer\'edi's theorem for arithmetic progressions of length four.
\newblock {\em Geom. Funct. Anal.}, 8(3):529--551,1998. 

\bibitem[Gr]{Gr}
A. Granville. 
\newblock An introduction to additive combinatorics.
\newblock Additive combinatorics, 1--27, CRM Proc. Lecture Notes, 43, Amer. Math. Soc., Providence, RI, 2007. 

 \bibitem[HaRu]{HaRu}
 Y. O. Hamidoune and J. Ru\'e. 
\newblock A lower bound for the size of a Minkowski sum of dilates.
\newblock {\em  Combin. Probab. Comput.}, 20(2):249--256, 2011.

\bibitem[HRY]{HRY}
F. Hennecart, G. Robert, and A. Yudin. 
\newblock On the number of sums and differences.
\newblock {\em  Structure theory of set addition. Ast\'erisque}, 258:173--178, 1999. 

\bibitem[KaKo]{KaKo}
N. H. Katz and P. Koester. 
\newblock On additive doubling and energy.
\newblock {\em  SIAM J. Discrete Math.}, 24(4):1684--1693, 2010.

\bibitem[KaSh]{KaSh}
N. H. Katz and C.-Y. Shen. 
\newblock A slight improvement to Garaev's sum product estimate.
\newblock {\em  Proc. Amer. Math. Soc.}, 136(7):2499--2504, 2008. 

\bibitem[Mu]{Mu}
A. Mudgal. 
\newblock Sums of linear transformations in higher dimensions.
\newblock {\em  Q. J. Math.}, 70(3):965--984, 2019. 

\bibitem[MuPe]{MuPe}
B. Murphy and G. Petridis. 
\newblock Products of differences over arbitrary finite fields.
\newblock {\em   Discrete Anal.}, Paper No. 18, 42 pp., 2019. 

\bibitem[Na]{Na}
M. B. Nathanson. 
\newblock Inverse problems for linear forms over finite sets of integers.
\newblock {\em  J. Ramanujan Math. Soc.} 23(2):151--165, 2008.

\bibitem[Pe]{Pe}
G. Petridis. 
\newblock New proofs of Pl\"unnecke-type estimates for product sets in groups.
\newblock {\em  Combinatorica}, 32(6):721--733, 2012. 

\bibitem[Pla]{Pla}
A. Plagne. 
\newblock Sums of dilates in groups of prime order.
\newblock {\em  Combin. Probab. Comput.}, 20(6):867--873, 2011.

\bibitem[PlTr]{PlTr}
A. Plagne and S. Tringali. 
\newblock Sums of dilates in ordered groups.
\newblock {\em  Comm. Algebra}, 44(12):5223--5236, 2016. 

\bibitem[Plu]{Plu} 
H. Pl\"unnecke. 
\newblock Eine zahlentheoretische Anwendung der Graphentheorie.
\newblock {\em J. Reine Angew. Math.}, 243:171--183, 1970. 

\bibitem[Ru1]{Ru1}
I. Z. Ruzsa. 
\newblock On the cardinality of ${A}+{A}$ and ${A}-{A}$.
\newblock  Combinatorics (Proc. Fifth Hungarian Colloq., Keszthely, 1976), Vol. II, pp. 933--938, Colloq. Math. Soc. J\'anos Bolyai, 18, North-Holland, Amsterdam-New York, 1978. 

\bibitem[Ru2]{Ru2}
I. Z. Ruzsa. 
\newblock An application of graph theory to additive number theory.
\newblock {\em Sci. Ser. A Math. Sci. (N.S.)}, 3:97--109, 1989. 


\bibitem[Ru3]{Ru3}
I. Z. Ruzsa. 
\newblock Cardinality questions about sumsets.
\newblock Additive combinatorics, 195--205, CRM Proc. Lecture Notes, 43, Amer. Math. Soc., Providence, RI, 2007. 


\bibitem[Ru4]{Ru4}
I. Z. Ruzsa. 
\newblock Sumsets and structure.
\newblock Combinatorial number theory and additive group theory, 87--210, Adv. Courses Math. CRM Barcelona, Birkh\"auser Verlag, Basel, 2009. 

\bibitem[Sch]{Sch}
T. Schoen. 
\newblock Near optimal bounds in Freiman's theorem.
\newblock {\em  Duke Math. J.}, 158(1):1--12, 2011.

\bibitem[ScSh]{ScSh}
T. Schoen and I. D. Shkredov. 
\newblock Higher moments of convolutions.
\newblock {\em  J. Number Theory}, 133(5):1693--1737, 2013.

\bibitem[Sha]{Sha}
G. Shakan. 
\newblock Sum of many dilates.
\newblock {\em  Combin. Probab. Comput.}, 25(3):460--469, 2016. 

\bibitem[Shk1]{Shk1}
I. D. Shkredov. 
\newblock Some new results on higher energies.
\newblock {\em Trans. Moscow Math. Soc.}, 2013, 31--63.

\bibitem[Shk2]{Shk2}
I. D. Shkredov. 
\newblock Some new inequalities in additive combinatorics.
\newblock {\em  Mosc. J. Comb. Number Theory}, 3(3-4):189--239, 2013. 

\bibitem[Shk3]{Shk3}
I. D. Shkredov. 
\newblock Energies and structure of additive sets.
\newblock {\em  Electron. J. Combin.}, 21, Paper 3.44, 53 pp., 2014. 

\bibitem[SSV]{SSV}
B. Sudakov, E. Szemer\'edi and V. H.  Vu. 
\newblock On a question of Erd\H{o}s and Moser.
\newblock {\em  Duke Math. J.}, 129(1):129--155, 2005. 

\bibitem[TaVu]{TaVu}
T. Tao and V. H . Vu. 
\newblock Additive combinatorics.
\newblock Studies in Advanced Mathematics, 105. Cambridge University Press, Cambridge, 2010. xviii+512 pp.

\end{thebibliography}


\begin{dajauthors}
\begin{authorinfo}[pgom]
  Brandon Hanson\\
  Department of  Mathematics \& Statistics, The University of Maine,\\ Orono, ME 04469, USA \\
  brandon\imagedot{}w\imagedot{hanson}\imageat{}gmail\imagedot{}com \\
  \url{http://www.brandonhanson.ca}
\end{authorinfo}
\begin{authorinfo}[johan]
  Giorgis Petridis\\
  Department of Mathematics, University of Georgia,\\
  Athens, GA 30602, USA\\
  giorgis\imageat{}cantab\imagedot{}net \\
  \url{https://faculty.franklin.uga.edu/petridis/}
\end{authorinfo}
\end{dajauthors}

\end{document}